\documentclass[twoside,11pt,reqno]{amsart}
\usepackage{macros}

\usepackage[capitalise]{cleveref}

\DeclarePairedDelimiterX{\set}[1]{\{}{\}}{\setargs{#1}}
\NewDocumentCommand{\setargs}{>{\SplitArgument{1}{;}}m}
{\setargsaux#1}
\NewDocumentCommand{\setargsaux}{mm}
{\IfNoValueTF{#2}{#1} {#1\,\delimsize|\,\mathopen{}#2}}

\newlist{thmenum}{enumerate}{1}
\setlist[thmenum]{label=(\roman*), ref=\thethm.(\roman*)}
\crefalias{thmenumi}{thm}

\newlist{propenum}{enumerate}{1}
\setlist[propenum]{label=(\roman*), ref=\theprop.(\roman*)}
\crefalias{propenumi}{prop}

\newlist{corenum}{enumerate}{1}
\setlist[corenum]{label=(\roman*), ref=\thecor.(\roman*)}
\crefalias{corenumi}{cor}

\newlist{lemenum}{enumerate}{1}
\setlist[lemenum]{label=(\roman*), ref=\thelm.(\roman*)}
\crefalias{lemenumi}{lm}

\newlist{exampleenum}{enumerate}{1}
\setlist[exampleenum]{label=(\alph*), ref=\theexamples.(\alph*)}
\crefalias{exampleenumi}{example}

\newlist{remarksenum}{enumerate}{1}
\setlist[remarksenum]{label=(\roman*), ref=\theremarks.(\roman*)}
\crefalias{remarksenumi}{remark}

\newlist{defenum}{enumerate}{1}
\setlist[defenum]{label=(\alph*), ref=\theDef.(\alph*)}
\crefalias{defenumi}{Def}

\begin{document}

\date{\today}

\title{Perfect loops of proper schemes}
\author{Alexander B. Ivanov}
\address{Fakult\"at f\"ur Mathematik, Ruhr-Universit\"at Bochum, Universit\"atsstrasse 150, D-44780 Bochum, Germany}
\email{a.ivanov@rub.de}

\maketitle

\begin{abstract}
We study the perfect loop space $LX$ of  a proper scheme $X$ over a local field $F$. We show that if $X$ contains a rational curve, then $LX$ is not ind-representable, even in a very weak sense. We also show that if $X$ admits a proper model $\fX$ over the integers $\caO$ of $F$, whose geometric special fiber does not contain rational curves, then $L^+\fX \stackrel{\sim}{\to} LX$ is an isomorphism via the natural map, where $L^+\fX$ denotes the perfect positive loop space (jet scheme) of $\fX$. In particular, $LX$ is a scheme in this case.
\end{abstract}

\section{Introduction}
Let $p$ be a prime.
Let $k$ be a perfect field of characteristic $p>0$ and let $F$ be either $k(\!(\varpi)\!)$ or a finite totally ramified extension of $W(k)[p^{-1}]$, where $W(\cdot)$ denotes the $p$-typical Witt-vectors. For example, $F$ could be any non-Archimedean local field.

If $X$ is an affine scheme of finite type over $F$, then its loop space $LX$ is a strict ind-scheme over $k$. In \cite[Theorem 5.1]{Ivanov_arc} it was shown that if $X$ is quasi-projective, then the restriction of $LX$ to perfect $k$-algebras is an arc-sheaf in the sense of \cite{BhattM_18}. After generalizing this to all separated $F$-schemes, in this note we study $LX$ for  schemes $X$ \emph{proper} over $F$. We regard $LX$ as a functor on perfect $k$-algebras. We show that there is an (incomplete) dichotomy: if $X$ contains a rational curve, $LX$ drastically fails to be ind-representable. On the other hand, if $X$ admits a proper model $\fX$ over the integers of $F$ without rational curves in the special fibre, then surprisingly the natural map $L^+\fX \stackrel{\sim}{\to} LX$ is an isomorphism, where $L^+\fX$ denotes the positive loop space (or jet scheme) of $\fX$; in particular, $LX$ is representable. Let us mention that loops of non-affine schemes were only rarely considered in the literature. A variant involving only nilpotent thickenings (hence somewhat orthogonal to our approach) was considered by Kapranov--Vasserot \cite{KapranovV_04}; more recently, $L(G/B)$ for $G$ a reductive group over $F$ and $B \subseteq G$ a Borel subgroup played a role in \cite{Ivanov_arc}.
We now describe our results in more detail.

\subsection*{Loop spaces are arc-sheaves.} The arc-topology on schemes was introduced in \cite{BhattM_18}. Among other things, it was shown there that when restricted to affine perfect schemes over $\bF_p$ the arc-topology is subcanonical. As mentioned above, \cite[Theorem 5.1]{Ivanov_arc} shows that $LX$ is an arc-sheaf when $X$ is a quasi-projective scheme over $F$. Generalizing this, we prove in \Cref{prop:loop_arc} that for a separated $F$-scheme $X$, $LX$ is an arc-sheaf. This follows by combining Tannakian reconstruction (of a qcqs scheme from the symmetric monoidal $\infty$-category of its perfect complexes) \cite[Theorem 1.5]{Bhatt_16} with two descent results: v-descent for perfect complexes  in a perfectoid setting \cite[Theorem 2.1]{AnschuetzLB_25_Fourier} (plus identification of perfect complexes on an affinoid perfectoid with perfect complexes over its global sections via \cite{Andreychev_21}) and a similar result in a more (derived-)schematic setting, \cite[Theorem 1.7]{Mathew_FaithfullyFlatDescent_22}.

\subsection*{Loops of the projective line}
Unsurprisingly, $L\bfP^1$ turns out to be not ind-representable, even when this term is understood in the very weak sense of being a filtered colimit of arbitrary (perfect) schemes along arbitrary transition maps (see \Cref{prop:P1_nonrep}). As a consequence (\Cref{cor:non_ind_rep_in_general}) we deduce that $LX$ is not ind-representable if there exists a non-constant map $\bfP^1 \to X$. In the equal characteristic case (including the case $\charac(k)=0$) Maffei gives a different proof of the non-ind-representability of $L\bfP^1$ in the forthcoming work \cite{Maffei_unpub}.

In the recent work \cite{MaffeiMV_26}, Maffei--Melani--Vezzosi investigate versions of the affine Grassmannian for a solvable group $G$ over a two-dimensional local field of equal characteristic. Theorem A of \emph{loc.cit.} states that all of them are ind-representable. From \Cref{prop:P1_nonrep} (as well as from \cite{Maffei_unpub}) it follows (\Cref{cor:MaffeiMV_indrep}) that for one of these versions -- the loop space ${}^L\mathcal{G}r_{G}=L(\mathcal{G}r_G)$ of the usual affine Grassmannian -- their result is essentially optimal: if $G$ is reductive, then ${}^L\mathcal{G}r_{G}$ is ind-representable if and only if $G$ is a torus.

\subsection*{Rational curves in the special fibre}

Quite surprisingly, it turns out that the behaviour described above is not typical for other proper schemes. Let $\caO$ denote the integers of $F$ and let $\fX$ be a proper scheme over $\caO$ with generic fibre $X/F$. One can form the (perfected) positive loop space $L^+\fX$ (also called jet scheme), and there is a natural map $L^+\fX \to LX$. Note that $L^+\fX$ is the cofiltered inverse limit of truncated jet schemes along affine transition maps, and therefore a rather explicit object; in particular, it is representable. For a detailed study of truncated jet schemes we refer to \cite{BertapelleG_18_arxiv}. Let $\bar k$ be a fixed algebraic closure of $k$.

\begin{introThm}[see \Cref{thm:Lplus_equals_L}]
Let $\fX$ be a proper $\caO$-scheme.
Suppose that any $\bar k$-morphism $\bfP^1_{\bar k} \to \fX_{\bar k}$ is constant. Then the natural map $L^+\fX \stackrel{\sim}{\to} LX$ is an isomorphism. In particular, $LX$ is representable by a perfect scheme.
\end{introThm}

We note that \cite[Theorem 1.1]{Anschuetz_22} is an analogue of this result for torsors under parahoric group schemes, see \Cref{rem:Johannes}. We deduce the remarkable consequence that the functor $L^+$ on the category $\caC_\caO$ of proper schemes satisfying the condition of the theorem factors through the generic fiber:
\begin{introCor}[see \Cref{cor:isogeny}]
The jet scheme functor \[\fX \mapsto L^+\fX \,\, \colon \,\,\caC_\caO \to \{\text{perfect $k$-schemes}\}\] factors through the generic fibre functor $\fX \mapsto \fX \times_\caO F$ from $\caC_\caO$ to the full subcategory of all proper $F$-schemes which admit a model $\fX \in \caC_\caO$.
\end{introCor}

\noindent \textbf{Outline of proof of \Cref{thm:Lplus_equals_L}.} Using that $L^+\fX$ and $LX$ are both arc-sheaves (by \cite[Theorem 5.17]{BhattM_18} and \Cref{prop:loop_arc} respectively), we reduce to showing that $L^+\fX(V) \to LX(V)$ is bijective for complete absolutely integrally closed valuation rings $V$ of rank one.

Let $\bW(V)$ be $V[\![\varpi]\!]$ if $\charac(F)=p$, resp. the ramified Witt-vectors of $V$ otherwise. Then $\Spec \bW(V)$ is a local scheme with closed point $s$. Essentially by Stein factorization the above bijectivity claim reduces to the fact that for a proper birational morphism $\Gamma \to \Spec\bW(V)$
which is an isomorphism outside $s$, any positive-dimensional irreducible component of the reduced fibre over $s$ is a rational curve. The generic point $\eta$ of such a component gives rise to a rank one point of the adic space $\caY = \Spa(\bW(V),\bW(V)) \sm \{[t] = 0\}$, where $t$ is a pseudo-uniformizer of $V$. This is exactly the adic space appearing in the Fargues--Fontaine construction, see \cite[Chapter II]{FarguesScholze}. The residue field at $\eta$ embeds into the secondary residue field at the corresponding point of $\caY$. In the equal characteristic case $\caY$ is just the open unit disc over the fraction field of $V$, its points are explictly  parameterized, the only points whose secondary residue fields are non-algebraic over the residue field $\kappa$ of $V$ are the points of type (2), in which case the residue field is isomorphic to $\kappa(z)$. We can conclude as any non-algebraic subextension of $\kappa(z)/\kappa$ is rational over $\kappa$. The mixed characteristic case reduces to the above by passing to a perfectoid cover of $\caY$ and tilting (as tilting preserves secondary residue fields).

\subsection*{Acknowledgements} The author wants to thank Ian Gleason and Jack Sempliner for helpful discussions on this paper, Ian Gleason for pointing out the observation recorded in Remark 3.2, and Ulrich G\"ortz and Peter Scholze for helpful suggestions how to improve the manuscript. The author also wants to thank Timo Richarz, Peter Scholze and Xinwen Zhu for discussions he had with them in the past, which motivated him to investigate the ind-representability properties of $L\bfP^1$.

\subsection*{Statement of AI use}
The author used extended discussions with generative AI while preparing this article.
Apart from searching references, helping to find proof details of technical lemmas, and reviewing arguments,
the models contributed the following substantial mathematical ideas:
\begin{itemize}
\item to use the $T_0$-property of the topological space underlying a scheme to derive a contradiction in the proof of \Cref{prop:P1_nonrep};
\item to use the Stein factorization and the closure of the graph in the proof of \Cref{thm:Lplus_equals_L} and use Temkin's results on (semi-)stable models (see \Cref{rem:Temkin_eq_char_case}): in particular, the original suggestion that $L^+\fX=LX$ should often hold is due to the model;
\item the final proof of \Cref{thm:Lplus_equals_L} was developed by the author jointly with the model (some of the ideas are due to the author, some are due to the model);
\item essentially the full proof of \Cref{prop:loop_arc} (to use Tannakian reconstruction via perfect complexes and their descent properties).
\item to apply \Cref{prop:P1_nonrep} to \cite{MaffeiMV_26} in \Cref{cor:MaffeiMV_indrep}.
\end{itemize}

\section{Notation and preliminaries}

We let $p$, $k$, $F$ and $\caO$ be as in the introduction.
Let $\Perf_{k}$ be the category of perfect $k$-algebras. We denote by $\varpi \in \caO$ a uniformizer of $F$. For $R \in \Perf_k$ let $W(R)$ be the $p$-typical Witt vectors of $R$. There is a unique $\varpi$-adically complete, separated $\caO_F$-algebra $\bW(R)$ which satisfies $\bW(R)/\varpi\bW(R) = R$, and in which $\varpi$ is a non-zero divisor. Explicitly,
\[
\bW(R) = \begin{cases} W(R) \otimes_{W(k)} \caO_F & \text{if ${\rm char}(F) = 0$} \\ R[\![\varpi]\!]& \text{if ${\rm char}(F) = p$.} \end{cases}
\]
In particular, $\bW(k)[\varpi^{-1}] = F$. We fix an algebraic closure $\bar{k}$ of $k$ and put $\breve F = \bW(\bar k)[\varpi^{-1}]$, the $\varpi$-adic completion of the maximal unramified extension of $F$.

For $R \in \Perf_k$ we will use the notations \[ \caY_R = \Spec \bW(R) \quad \text{ and }\quad \widetilde R = \bW(R)[\varpi^{-1}].\]

Finally, note that for $R \in \Perf_k$ any element $x \in \widetilde R$ admits a unique expansion $x = \sum_{i=-N}^{\infty} [x_i]\varpi^i$ with $x_i \in R$, where $[x_i] = x_i$ if ${\rm char}(F) > 0$ and $[x_i] \in \bW(R)$ is the Teichm\"uller lift if ${\rm char}(F) = 0$. Conversely, any such series defines an element of $\widetilde R$. Also $x \in \bW(R)$ if and only if $x_i = 0$ for $i<0$; for $x \in \bW(R)$, $x \in \bW(R)^\times \Leftrightarrow x_0 \in R^\times$.

\subsection{Ind-schemes.} In this note, by a \emph{weak ind-(perfect scheme)} (or simply \emph{weak ind-scheme}) we mean a set-valued functor of $\Perf_k$ which can be written as a termwise filtered colimit of perfect schemes. We do not require these schemes to be quasi-compact and we also do not pose any conditions on the transition maps. We refer to \cite{BhattM_18} for the definition of the arc-topology on $\Perf_k$.

\begin{rem}\label{rem:filtered_colims}
\begin{itemize}
\item[(i)] If $(F_i)_i$ is a filtered system of arc-sheaves on $\Perf_k$, then the termwise colimit $F=\dirlim_i F_i$ formed in presheaves, is again an arc-sheaf. Indeed, the arc-topology on $\Perf_k$ is finitary by \cite[Remark 1.4]{BhattM_18}, hence any arc-cover in $\Perf_k$ can be refined by a singleton affine cover. Let $\Spec B \to \Spec A$ be such a cover. Then $F(A)= \dirlim_i {\rm Eq}(F_i(B) \rightrightarrows F_i(B\otimes_A B)) \simeq {\rm Eq}(F(B) \rightrightarrows F(B\otimes_A B))$, as filtered colimits commute with finite limits in sets. The same holds for the fpqc-topology.
\item[(ii)] By \cite[Theorem 5.17]{BhattM_18} any perfect qcqs $k$-scheme regarded as a presheaf on $\Perf_k$ via Yoneda is an arc-sheaf. By Gabber's observation \cite[Remark 4.6]{Bhatt_16} this generalizes to arbitrary perfect $k$-schemes $X$. Indeed, we may write $X(A) = \dirlim_{U\subseteq X \text{qc open}} U(A)$ for any $A\in \Perf_k$. From this and part (i) of the remark, we are reduced to the case that $X$ is quasi-compact. But then by \cite[Tag 03K0]{StacksProject} (see also \cite[Remark 4.6]{Bhatt_16}) we may write $X=\colim_i X_i$ as a filtered colimit in fpqc-sheaves on $\Perf_k$, where the $X_i$ are qcqs perfect schemes. By part (i) of the remark this colimit can be computed termwise and agrees with the colimit in arc-sheaves. Thus $X$ is an arc-sheaf.
\end{itemize}
\end{rem}

\begin{lm}\label{lm:weak_ind}
Let $X=\dirlim_i X_i$ be a weak ind-scheme. Then
\begin{itemize}
\item[(i)] If $T$ is a qcqs perfect $k$-scheme (via Yoneda regarded as a functor on $\Perf_k$), then $\dirlim_i\Hom(T,X_i) = \Hom(T,X)$.
\item[(ii)] $X$ is an arc-sheaf on $\Perf_k$.
\end{itemize}
\end{lm}
\begin{proof}
We leave part (i) as an exercise. Part (ii) follows from \Cref{rem:filtered_colims}(i). \qedhere
\end{proof}

Note that part (ii) of the lemma is very similar to \cite[Lemma 1.4]{Richarz_20_basics}. There the stronger condition that $X_i \to X$ is a monomorphism (which is not available in our setting) is only used to verify the sheaf condition on non-quasi-separated test objects (cf. Remark 1.5 of \emph{loc.cit.}).

We call a functor $X$ on $\Perf_k$ a \emph{strict ind-scheme} if it can be written as a termwise filtered colimit $X = \dirlim_i X_i$ along closed immersions with all $X_i$ qcqs perfect schemes (this is the perfected version of the ind-scheme notion of \cite[\S6.3.1]{Drinfeld_03}). Clearly, any strict ind-scheme is a weak ind-scheme.

We say a functor on $\Perf_k$ is \emph{representable} if it is representable by a perfect $k$-scheme; we say it is \emph{(weakly resp. strictly) ind-representable} if it is representable by a (weak resp. strict) ind-scheme.

\subsection{Loop spaces} For a scheme $X$ over $F$ the (perfect) \emph{loop functor} of $X$ is defined by
\[
LX \colon \Perf_k \to {\rm Sets}, \quad R \mapsto X(\bW(R)[\varpi^{-1}]).
\]
If $X/F$ is affine and of finite type over $F$, then $LX$ is well-known to be strictly ind-representable, see \cite[\S1.a]{PappasR_08} and \cite[Proposition 1.1]{Zhu_17}. If $X$ is quasi-projective over $F$, then \cite[Theorem 5.1]{Ivanov_arc} guarantees at least that $LX$ is an arc-sheaf on $\Perf_k$.

Similarly, for a scheme $\fX$ over $\caO$ the \emph{positive loop functor} (or \emph{jet-scheme}) of $\fX$ is defined by
\[
L^+\fX \colon \Perf_k \to {\rm Sets}, \quad R \mapsto \fX(\bW(R)).
\]
For $r\geq 1$, there is also its truncated version $L^+_r\fX(R) = \fX(\bW(R)/\varpi^r)$. Then $L^+\fX = \prolim_r L^+_r\fX$, along affine transition maps: this follows from \cite[Theorem 4.1 and Remarks 4.3, 4.6]{Bhatt_16} and \cite[Proposition 9.23]{BertapelleG_18_arxiv}. In particular, $L^+\fX$ is a (perfect) scheme over $k$, see \cite{BertapelleG_18_arxiv}.

If $\fX$ is an $\caO$-scheme and $X = \fX \times_{\caO} F$ is its generic fibre, then restriction along $\Spec \bW(R)[\varpi^{-1}] \to \Spec \bW(R)$ defines a natural morphism
\[
\iota = \iota_{\fX} \colon L^+\fX \to LX.
\]
Let us record the following.
\
\begin{lm}\label{lm:injective}
If $\fX$ is separated over $\caO$, then $\iota$ is injective.
\end{lm}
\begin{proof}
Given two points $f,g : \caY_R \to \fX$ in $L^+\fX(R)$, their equalizer ${\rm Eq}_{f,g} \subseteq \caY_R$ equals the pullback of $(f,g) : \caY_R \to \fX \times_{\caO} \fX$ along the diagonal of $\fX \to \Spec \caO$, hence a closed subscheme of $\caY_R$; now $\iota(f) = \iota(g)$ just means that ${\rm Eq}_{f,g}$ contains $\Spec \widetilde R$, and hence equals $\caY_R$.
\end{proof}

\subsection{Valuation rings of rank one}

Valuation rings of rank $\leq 1$ behave differently from higher rank valuation rings in several respects. We record two facts, both failing in higher ranks.

\begin{lm}\label{lm:normality}
Let $V \in \Perf_k$ be a valuation ring of rank $\leq 1$. Then $\bW(V)$ is a normal domain.
\end{lm}
\begin{proof}
Such $V$ is completely integrally closed (this is immediate). If ${\rm char}(F) > 0$, then by \cite[Exercise 9.5 on p.69]{Matsumura_CRT}, $\bW(V) = V[\![\varpi]\!]$ is also completely integrally closed, and hence an integrally closed domain. If ${\rm char}(F) = 0$, then $\bW(V)$ is an integrally closed domain by \cite[Proposition 5.4(iii)]{Shimomoto_17}. Note that \cite[Proposition 5.4]{Shimomoto_17} assumes finiteness of $k$; however, the proof of part (iii) (and the reference \cite[Corollary 5.10]{Shimomoto_15} used in the proof) of this proposition does not use this assumption.\qedhere
\end{proof}

The above lemma fails for valuation rings of rank $>1$ (see \cite[Exercise 10.4 on p. 77]{Matsumura_CRT} and references therein). On the other hand, we crucially make use of it in the proof of \Cref{thm:Lplus_equals_L} below; thus it is really important to work with arc-sheaves and not with v-sheaves.

We also record the following flatness result, which we could not find anywhere in the literature. In contrast to \Cref{lm:normality} we do not use it in the proof of \Cref{thm:Lplus_equals_L}. In the equal characteristic case, there is a purely schematic proof of this theorem through Temkin's semistable modifications \cite{Temkin_sst}, in which case this flatness simplifies the proof; see \Cref{rem:Temkin_eq_char_case}. Note that this lemma also fails for valuation rings of rank $>1$, see \cite[0F1Z]{StacksProject}.

\begin{lm}\label{lm:flatness_rk1val}
Let $V$ be a valuation ring of rank $1$. Then $V[\![\varpi]\!]$ is flat over $V[\varpi]$.
\end{lm}
\begin{proof}
Let $\fm \subseteq V$ denote the maximal ideal and let $\fm[\varpi] \subseteq V[\varpi]$ denote the ideal generated by $\fm$. By \cite[0H7N]{StacksProject} applied to the non-zero divisor $\varpi$ it suffices to prove that $V(\!(\varpi)\!)$ is flat over $V[\varpi^{\pm 1}]$, or equivalently, that $V(\!(\varpi)\!)$ is flat over the localization $(V[\varpi] \sm (\fm,\varpi))^{-1} V[\varpi^{\pm 1}] = V[\varpi]_{(\fm,\varpi)}[\varpi^{-1}] =: A$. We claim that $A$ is a Pr\"ufer ring. Note that $\Spec A \subseteq \Spec V[\varpi]$ is a pro-(qc open) subset, any prime $\fq$ of $A$ determines a prime $\fp = V[\varpi] \cap \fq$ of $V[\varpi]$ satisfying $\fp \subseteq (\fm,\varpi)$ and $\varpi \not\in \fp$, and we have $V[\varpi]_\fp = A_{\fq}$. Put $K={\rm Frac}(V)$. If $\fp \cap V = \{0\}$, then $V[\varpi]_\fp$ is a localization of $K[\varpi]$ and hence a valuation ring. Otherwise $\fp \cap V = \fm$, so that $\fm[\varpi] \subseteq \fp \subseteq (\fm,\varpi)$. Reducing modulo $\fm$, $\fp/\fm[\varpi] \subseteq (\varpi) \subseteq (V/\fm)[\varpi]$. But $\varpi \not\in \fp$, so that we deduce $\fp = \fm[\varpi]$. But $V[\varpi]_{\fm[\varpi]}$ is the valuation ring of the Gaussian valuation (see e.g. \cite[Chap. VI, \S10]{Bourbaki_CommAlg}).
Thus all local rings of $A$ are valuation rings and hence $A$ is Pr\"ufer. Finally, $A \to V(\!(\varpi)\!)$ is injective ($A$ is a localization of $V[\varpi]$, hence contained in ${\rm Frac}(V[\varpi]) \subseteq K(\!(\varpi)\!)$, and $V(\!(\varpi)\!) \subseteq K(\!(\varpi)\!)$ too), and as $V(\!(\varpi)\!)$ is a domain it follows that $V(\!(\varpi)\!)$ is torsion-free over $A$. As $A$ is Pr\"ufer, this implies that $V(\!(\varpi)\!)$ is $A$-flat.
\end{proof}

\section{Loops of separated schemes are arc-sheaves}

The goal of this section is to prove that $LX$ is an arc-sheaf when $X$ is a separated scheme over $F$. We begin with some lemmas. Let ${\rm Perfd}_F$ denote the category of affinoid perfectoid spaces over $F$. Then ${\rm Perfd}_F$ is equipped with the arc-topology, which satisfies the usual finitarity property and where a map $\Spa(B,B^+) \to \Spa(A,A^+)$ is a cover if it is surjective on rank $1$ points; see \cite[Definition 3.1 and discussion after Corollary 12.2]{Scholze_Berkovich_motives} and \cite[Definition 5.6]{Ivanov_arc}. For a ring $R$, let ${{\caP}erf}(R)$ denote the $\infty$-category of  perfect complexes of $R$-modules. For a perfectoid space $Y$ let ${{{\caP}erf}}(Y)$ denote the $\infty$-category of perfect complexes of $\caO_Y$-modules on $Y$ (i.e., those which are locally for the analytic topology on $Y$ quasi-isomorphic to a bounded complex of finite locally free $\caO_Y$-modules).

\begin{lm}\label{lm:algebraiz}
Let $X$ be a quasi-separated scheme over $F$. The functor $X^{\Diamond_{\rm pre}} : \Spa(A,A^+) \mapsto X(A)$ is an arc-sheaf on ${\rm Perfd}_F$.
\end{lm}
\begin{proof}
Write $X=\dirlim_{i} X_i$ as a filtered colimit with $X_i$ qcqs open subsets. As $\Spec A$ is quasi-compact, we have
\[
X^{\Diamond_{\rm pre}}(\Spa(A,A^+)) = X(A) = \dirlim_i X_i(A) = \dirlim_i X_i^{\Diamond_{\rm pre}}(\Spa(A,A^+)) = (\dirlim_i X_i^{\Diamond_{\rm pre}})(\Spa(A,A^+)),
\]
where the last colimit is formed in the category of set-valued presheaves on ${\rm Perfd}_F$. Once we have shown that all $X_i^{\Diamond_{\rm pre}}$ are arc-sheaves, this colimit agrees with the colimit formed in the category of arc-sheaves on ${\rm Perfd}_F$ (any arc-covering in ${\rm Perfd}_F$ admits a refinement by singleton affinoid covering; hence filtered colimits in arc-sheaves can be computed termwise and agree with colimits in presheaves).

Thus we may assume that $X$ is qcqs over $F$. Moreover, it suffices to show that $X^{\Diamond_{\rm pre}}$ is a v-sheaf; indeed, then the independence of $X^{\Diamond_{\rm pre}}(\Spa(A,A^+))$ on $A^+$ immediately gives that $X^{\Diamond_{\rm pre}}$ is partially proper and then \cite[Corollary 12.2 and the discussion after it]{Scholze_Berkovich_motives} implies that $X^{\Diamond_{\rm pre}}$ is an arc-sheaf. Let $\Spa(B,B^+) \to \Spa(A,A^+)$ be a v-cover in ${\rm Perfd}_F$ and let $B_n$ denote the global sections of the $n$th term of the \v{C}ech-nerve of this covering (which is again affinoid perfectoid). By \cite[Theorem 5.43, Lemmas 5.46, 5.47]{Andreychev_21}, the natural functor ${{\caP}erf}(A) \to {{\caP}erf}(\Spa(A,A^+))$ is an equivalence of categories (see also \cite[p. 3655]{AnschuetzLB_25_Fourier}). Using
\cite[Theorem 2.1]{AnschuetzLB_25_Fourier} for the equivalence in the middle, we deduce
\begin{equation}\label{eq:perf_A_descent}
{{\caP}erf}(A) \simeq {{\caP}erf}(\Spa(A,A^+)) \simeq \prolim_{[n] \in \Delta} {{\caP}erf}(\Spa(B_n,B_n^+)) \simeq \prolim_{[n] \in \Delta} {{\caP}erf}(B_n)
\end{equation}
If $Z,T$ are two schemes, let ${\rm Fun}^{\otimes,{\rm ex}}({\caP}erf(Z),{\caP}erf(T))$ denote the $\infty$-category of exact symmetric monoidal functors ${\caP}erf(Z) \to {\caP}erf(T)$. Moreover, if $Z,T$ are schemes over $F$, let ${\rm Fun}_{F\text{-sch}}^{\otimes,{\rm ex}}({\caP}erf(Z),{\caP}erf(T))$ denote the $\infty$-category of those exact symmetric monoidal functors ${\caP}erf(Z) \to {\caP}erf(T)$ which commute with the natural functors from ${\caP}erf(\Spec F)$. If $Z$ is qcqs, \cite[Theorem 1.5]{Bhatt_16} implies that pullback induces natural equivalences
\[
\Hom_{F}(T,Z) \simeq {\rm Fun}_{F\text{-sch}}^{\otimes,{\rm ex}}({\caP}erf(Z),{\caP}erf(T))
\]
(a posteriori it follows that ${\rm Fun}_{F\text{-sch}}^{\otimes,{\rm ex}}({\caP}erf(Z),{\caP}erf(T))$ is equivalent to a full subcategory of ${\rm Fun}^{\otimes,{\rm ex}}({\caP}erf(Z),{\caP}erf(T))$).
Using this for the first and the last equivalences, we compute:
\begin{align*}
X(A) &\simeq {\rm Fun}_{F\text{-sch}}^{\otimes,{\rm ex}}({\caP}erf(X),{\caP}erf(A)) \\
&\simeq {\rm Fun}_{F\text{-sch}}^{\otimes,{\rm ex}}({\caP}erf(X),\prolim_{[n]\in\Delta}{\caP}erf(B_n)) \\
&\simeq \prolim_{[n]\in \Delta} {\rm Fun}_{F\text{-sch}}^{\otimes,{\rm ex}}({\caP}erf(X),{\caP}erf(B_n))\\
&\simeq \prolim_{[n]\in \Delta} X(B_n).
\end{align*}
using \eqref{eq:perf_A_descent} for the second equivalence. This precisely says that $X^{\Diamond_{\rm pre}}$ satisfies v-descent.
It only remains to justify that we may pull out the limit from ${\rm Fun}_{F\text{-sch}}^{\otimes,{\rm ex}}({\caP}erf(X),-)$ in the third equivalence above. This follows for the category of all functors, as ${\rm Fun}({\caP}erf(X),-)$ is a right adjoint and hence commutes with all limits (see \cite[Remark 7.5.2.3]{Lurie_Kerodon}); by \cite[Theorem 1.1.4.4]{Lurie_HA}, the equivalence ${\rm Fun}({\caP}erf(X),\prolim_{[n]\in \Delta}{\caP}erf(B_n)) \stackrel{\sim}{\to} \prolim_{[n]\in \Delta} {\rm Fun}({\caP}erf(X),{\caP}erf(B_n))$ restricts to an equivalence of the full subcategories of exact functors. It preserves monoidal structures by \cite[Corollary 3.2.2.5]{Lurie_HA}.
\end{proof}

Note that for affine schemes $X$ \Cref{lm:algebraiz} follows directly from \cite[Theorem 8.7]{Scholze_ECD} (and \cite{Scholze_Berkovich_motives}).

\begin{rem}
As was pointed out to the author by I. Gleason, \Cref{lm:algebraiz} fails in general for the small diamond functor $X^{\diamond_{\rm pre}} \colon \Spa(A,A^+) \mapsto X(A^+)$ for non-affine $X$. In fact, $(\bfP^1_{\bF_p})^{ \diamond_{\rm pre}}$ is not even an analytic sheaf.
\end{rem}

\begin{lm}\label{lm:fpqc_varpi-descent}
Let $X$ be a quasi-separated scheme over $F$. Let $B \to B'$ be a map of $\varpi$-torsion-free, $\varpi$-adically complete $\caO$-algebras such that $B/\varpi \to B'/\varpi$ is faithfully flat. Then \[ X(B[\varpi^{-1}]) = {\rm Eq}(X(B'[\varpi^{-1}]) \rightrightarrows X((B' \widehat\otimes_B B')[\varpi^{-1}])).\]
\end{lm}
\begin{proof}
As in the proof of \Cref{lm:algebraiz}, we may reduce to the case that $X$ is qcqs. Let $B_n' = \left(\bigotimes_{B}^{n+1} B'\right)^\wedge_\varpi$ be the (non-derived) $\varpi$-adically completed tensor product. We claim that this agrees with the derived completed tensor product. We use induction on $n$ and apply \cite[Remark 2 after Theorem A.25]{GrosSQ_23} to the bounded $B$-algebra $B'$ (taking $f=0$, $g=\varpi$), and the derived completed tensor product $M=\widehat{\bigotimes}_{B}^{\bL, n} B'$ of $B$-modules. To be able to apply \emph{loc.cit.} we must guarantee that the complex $M$ is a discrete, completely flat $B$-module. Note that discreteness of $M$ is guaranteed by the induction hypothesis for $n\geq 1$ (and is obvious if $n=0$). For complete flatness (see Definition A.3 of \emph{loc.cit.}), note that
\[
M \otimes_B^\bL B/\varpi \simeq \left({\bigotimes}_{B}^{\bL, n} B'\right) \otimes_B^\bL B/\varpi \simeq {\bigotimes}_{B/\varpi}^{\bL, n} (B' \otimes_B^\bL B/\varpi ) \simeq {\bigotimes}_{B/\varpi}^{n} B'/\varpi,
\]
where the first equivalence follows by Proposition A.10(1) of \emph{loc.cit.}, and the last equivalence follows from that $B$ and $B'$ are $\varpi$-torsion free (compute the derived tensor product using the free resolution $B/\varpi \simeq [B \stackrel{\varpi}{\to} B]$ in degrees $-1,0$). Thus  the left hand side is discrete and a flat $B/\varpi$-module (by assumption of the lemma). This proves the claim.

By the above claim along with \cite[Theorem 1.7]{Mathew_FaithfullyFlatDescent_22} (which works with derived complete tensor products), we have that ${\caP}erf(B[\varpi^{-1}]) = \prolim_{[n] \in \Delta} {\caP}erf(B_n'[\varpi^{-1}])$. As in the previous lemma, \cite[Theorem 1.5]{Bhatt_16} shows
\[
X(B[\varpi^{-1}]) = {\rm Fun}_{F\text{-sch}}^{\otimes,{\rm ex}}({\caP}erf(X),{\caP}erf(B[\varpi^{-1}]))
\]
and similarly for the (ordinary) rings $B'_n[\varpi^{-1}]$. Now the same computation as in the end of the proof of \Cref{lm:algebraiz} (with $B[\varpi^{-1}],B_n'[\varpi^{-1}]$ instead of $A,B_n$) finishes the proof.
\end{proof}

We are ready to prove the main result of this section, which generalizes \cite[Theorem 5.1]{Ivanov_arc}.

\begin{prop}\label{prop:loop_arc}
Let $X$ be a separated scheme over $F$. Then $LX$ is an arc-sheaf on $\Perf_k$.
\end{prop}
\begin{proof}
Let $F_\infty$ be the completion of $F(\varpi^{\frac{1}{p^{\infty}}})$, with integers $\caO_{\infty}$. Then $F_\infty$ is a perfectoid extension of $F$. For $R\in \Perf_k$, write
\[\begin{alignedat}{3}
A_R^+ &= \bW(R) &\quad &\text{ and } \quad& A_R &= A_{R}^+[\varpi^{-1}] \\
A_{\infty,R}^+ &= \bW(R) \, \widehat\otimes_{\caO} \,\caO_\infty &\quad &\text{ and } \quad & A_{\infty,R} &= A_{\infty,R}^+[\varpi^{-1}]
\end{alignedat}\]
We equip $A_R^+$ and $A_{\infty,R}^+$ with $\varpi$-adic topology and regard them as open subrings of $A_R$, $A_{\infty,R}$ respectively. Then $(A_R,A_R^+)$ and $(A_{\infty,R},A_{\infty,R}^+)$ are Tate--Huber pairs, and we have the affinoid adic spaces
\[
U_R^\infty = \Spa(A_{\infty,R},A_{\infty,R}^+) = U_R \times_{\Spa(F,\caO)} \Spa(F_\infty,\caO_\infty) \longrightarrow U_R = \Spa(A_R,A_R^+),
\]
where $U_R^\infty$ is an affinoid perfectoid cover of $U_R$. This setting is the same as in \cite[proof of Proposition 5.9]{Ivanov_arc} (see also \cite[Proposition 19.5.3]{ScholzeW_20}).

Let $R \to R'$ be an arc-cover in $\Perf_k$, and let $R'' = R' \otimes_R R'$ (which is again in $\Perf_k$). Consider the following commutative diagram:
\[
\begin{tikzcd}[column sep=large, row sep=large]
X(A_R)
  \arrow[r] \arrow[d]
&
X(A_{R'})
  \arrow[r, shift left=0.5ex]
  \arrow[r, shift right=0.5ex]
  \arrow[d]
&
X(A_{R''})
  \arrow[d]
\\
X(A_{\infty,R})
  \arrow[r]
  \arrow[d, shift left=0.5ex]
  \arrow[d, shift right=0.5ex]
&
X(A_{\infty,R'})
  \arrow[r, shift left=0.5ex]
  \arrow[r, shift right=0.5ex]
  \arrow[d, shift left=0.5ex]
  \arrow[d, shift right=0.5ex]
&
X(A_{\infty,R''})
\\
X\!\left(
  \left(
    A_{\infty,R}^{+}
    \widehat{\otimes}_{A_R^{+}}
    A_{\infty,R}^{+}
  \right)[\varpi^{-1}]
\right)
  \arrow[r]
&
X\!\left(
  \left(
    A_{\infty,R'}^{+}
    \widehat{\otimes}_{A_{R'}^{+}}
    A_{\infty,R'}^{+}
  \right)[\varpi^{-1}]
\right)
&
{}
\end{tikzcd}
\]
Let $x' \in {\rm Eq}(X(A_{R'}) \rightrightarrows X(A_{R''}))$. We need to show that it comes from a unique element of $X(A_R)$. The image $x'_{\infty} \in X(A_{\infty,R'})$ of $x'$ lies in the equalizer of the two maps $X(A_{\infty,R'}) \to X(A_{\infty,R''})$. Exploiting \Cref{lm:derived_perfd_tensor} to identify $U_{R''}^\infty$ with $U_{R'}^\infty \times_{U_R^\infty} U_{R'}^\infty$, the middle row can be identified with
\[ X^{\Diamond_{\rm pre}}(U_R^{\infty}) \to X^{\Diamond_{\rm pre}}(U_{R'}^\infty) \rightrightarrows  X^{\Diamond_{\rm pre}}(U_{R'}^\infty \times_{U_R^\infty} U_{R'}^\infty).\]
Moreover, \cite[Proposition 5.3]{Ivanov_arc} implies that $U_{R'}^\infty \to U_R^\infty$ is an arc-cover. Thus, by \Cref{lm:algebraiz}, $x'_\infty$ comes from a unique element  $x_\infty \in X(A_{\infty,R})$. Moreover, the lower horizontal map is injective: indeed, as $X$ is a separated scheme, it suffices to show (by looking at the equalizer subscheme of two morphisms) that the ring map
\begin{align*}
(\bW(R) \widehat\otimes_\caO \caO_\infty \widehat\otimes_\caO \caO_\infty)[\varpi^{-1}] = (A_{\infty,R}^+ &\widehat{\otimes}_{A_{R}^+} A_{\infty,R}^+)[\varpi^{-1}] \to \\
&\to (A_{\infty,R'}^+ \widehat{\otimes}_{A_{R'}^+} A_{\infty,R'}^+)[\varpi^{-1}] = (\bW(R') \widehat\otimes_\caO \caO_\infty \widehat\otimes_\caO \caO_\infty)[\varpi^{-1}]
\end{align*}
is injective. Being an arc-cover in $\Perf_k$, $R \to R'$ is injective. Thus $\bW(R)/\varpi^n \to \bW(R')/\varpi^n$ is injective for all $n$. Tensoring twice with the flat map $\caO/\varpi^n \to \caO_\infty/\varpi^n$, passing to inverse limit over $n$ and inverting $\varpi$, leaves the resulting map injective.

It follows that the images of $x_\infty$ under the two vertical maps coincide, as this holds for $x'_\infty$ (as $x'_{\infty}$ comes from $X(A_{R'})$). Note that \Cref{lm:fpqc_varpi-descent} applies to $A_{R}^+ \to A_{\infty,R}^+$ as $A_{\infty,R}^+/\varpi = R \otimes_k \caO_\infty/\varpi$ is faithfully flat over $R$. Thus, $x_\infty$ comes from a unique element $x \in X(A_R)$. Moreover, the image of $x$ in $X(A_{R'})$ agrees with $x'$, as  the images of both in $X(A_{\infty,R'})$ coincide, and the middle column is exact (again by \Cref{lm:fpqc_varpi-descent}).
\end{proof}

\begin{lm}\label{lm:derived_perfd_tensor}
Let $R_1 \leftarrow R \to R_2$ be maps in $\Perf_k$. Then we have
\[ U^\infty_{R_1 \otimes_R R_2} \cong U^\infty_{R_1} \times_{U^\infty_R} U^\infty_{R_2}, \]
where the fiber product is formed either in adic or in perfectoid spaces.
\end{lm}

\begin{proof}
Comparing universal properties, the natural map $\bW(R_1) \widehat\otimes_{\bW(R)} \bW(R_2) \to \bW(R_1 \otimes_R R_2)$ is an isomorphism. From this, also $A_{\infty,R_1}^+
\widehat\otimes_{A_{\infty,R}^+} A_{\infty,R_2}^+ \to A_{\infty, R_1\otimes_R R_2}^+$ is an isomorphism (notation as in the proof of \Cref{prop:loop_arc}). Then the same holds after inverting $\varpi$, and the lemma follows.
\end{proof}

\section{$L\bfP^1$ is not weakly ind-representable}
Now suppose $X$ is proper over $F$ and let $\fX/\caO_F$ be a proper model. For any field $\ff\in \Perf_k$, $\bW(\ff)$ is a discrete valuation ring with fraction field $\widetilde \ff = \bW(\ff)[\varpi^{-1}]$, so the valuative criterion for properness gives
\begin{equation}\label{eq:1}
LX(\ff) = X(\,\widetilde \ff\,) = \fX(\bW(\ff)) = L^+\fX(\ff).
\end{equation}
That is, $\iota$ is bijective on field valued points.
However, for a valuation ring $V$ of rank $\geq 1$, one has $L\bfP^1(V) \neq L^+\bfP^1(V)$ (see proof of \Cref{prop:P1_nonrep}). Thus, heuristically, one should expect $L\bfP^1$ to be not representable by a scheme: otherwise  $L^+\bfP^1 \to L\bfP^1$ would be a strange map of perfect schemes which is a bijection on geometric points, but not an isomorphism. We confirm this heuristic:
\begin{prop}\label{prop:P1_nonrep}
$L\bfP^1$ is not weakly ind-representable.
\end{prop}

\begin{proof}
Suppose $L\bfP^1$ is weakly ind-representable and let $L\bfP^1 = \dirlim_{i\in I} X_i$ be an ind-presentation with $X_i$ arbitrary schemes. As $L^+\bfP^1$ is qcqs, there exists by \Cref{lm:weak_ind} some $n \in I$ such that $L^+\bfP^1 \to L\bfP^1$ factors through a morphism $\iota_i : L^+\bfP^1\to X_i$ for all $i\geq n$. As the composition $L^+\bfP^1 \to X_i \to L\bfP^1$ is bijective on all field valued points by \eqref{eq:1}, $L^+\bfP^1 \to X_i$ is injective on such points, and therefore (universally) injective by \cite[Tag 03MU]{StacksProject}. In particular, $|L^+\bfP^1| \to |X_i|$ is injective for all $i\geq n$, where $|\cdot|$ denotes the topological space of the underlying scheme.
Note that $L^+\bfP^1$ is irreducible (it is an iterated $\bA^1$-bundle over $\bP_k^1$), and let $\eta$ denote its generic point. It specializes to any other point of $L^+\bfP^1$, and this specialization relation remains true after applying $\iota_i$, that is for any $i\geq n$ and any $z \in |L^+\bfP^1|$, we have the specialization relation $\iota_i(\eta) \rightsquigarrow \iota_i(z)$.

Let $V \in \Perf_k$ be a rank one valuation ring and $f:\Spec V = \{\xi \rightsquigarrow s\} \to L\bfP^1$ a map. As $\Spec V$ is qcqs, $f$ factors through a map $f_m : \Spec V \to X_m$ for some $m\geq n$. Moreover, the restrictions of $f$ to both points of $\Spec V$ are field-valued points of $L\bfP^1$, hence factor through $L^+\bfP^1$ by \eqref{eq:1}. This gives two points $f(\xi), f(s) \in |L^+\bfP^1|$ and $f_m$ induces -- after enlarging $m$ if necessary -- a specialization relation between their images $\iota_m(f(\xi)) \rightsquigarrow \iota_m(f(s))$ in $X_m$. Suppose that we can construct $V$ and $f$, such that $\eta = f(s) \neq f(\xi)$. Then we would get specialization relations $\iota_m(f(\xi)) \rightsquigarrow \iota_m(f(s)) \rightsquigarrow \iota_m(f(\xi))$ in $|X_m|$ (first induced by $f_m$, second coming from $L^+\bfP^1$), but by injectivity of $|L^+\bfP^1| \to |X_m|$, $\iota_m(f(\xi)) \neq \iota_m(f(s))$. This contradicts the fact that $|X_m|$ is a $T_0$-space.

It remains to construct $V,f$ as above. We have the natural map $L^+\bfP^1 \to L^+_1 \bfP^1 = (\bfP_k^1)^{\rm perf}$. Denoting the homogeneous coordinates on $\bfP_k^1$ by $[x_0:y_0]$, let $D(x_0), D(y_0) \subseteq L^+\bfP^1$ be the preimages of the principal open subsets $\{x_0 \neq 0\}, \{y_0 \neq 0\} \subseteq \bfP^1$. Note that $D(x_0), D(y_0)$ are both affine and isomorphic to $L^+\bA^1 (\cong (\bA^\infty_k)^{\rm perf})$. In particular, $[x:y] \mapsto t:=\frac{y}{x} \colon D(x_0) \stackrel{\sim}{\to} L^+\bA^1$ is an isomorphism. Let $Z \subseteq D(x_0)$ be the closed subscheme given by $t_0 = 0$, where the unique expansion $t=\sum_{i=0}^{\infty}[t_i]\varpi^i$ determines the coordinates on $L^+\bA^1$.

Let $\kappa := k(a_0,a_1,a_2,\dots)^{\rm perf}$ with $a_i$ algebraically independent variables. Let $V = \kappa[\![t]\!]^{\rm perf}$ be a perfect rank one valuation ring with residue field $\kappa$ and fraction field $K={\rm Frac}(V)$. Write $\Spec V = \{\xi\rightsquigarrow s\}$. Let \[ f:= [t]\varpi^{-1} + [a_0] + [a_1]\varpi + [a_2] \varpi^2 + \dots \in \widetilde V. \]
Then $[f:1] \in \bfP^1(\widetilde V) = L\bfP^1(V)$. At the special point of $\Spec V$ one has $t=0$. As the $a_i$ are algebraically independent, the image of $[f(\kappa):1] : \Spec \kappa \to D(y_0) \subseteq L^+\bfP^1$ is exactly the generic point of $D(y_0)$ (the map $H^0(D(y_0),\caO) \to \kappa$ is injective),
which equals $\eta$ as $D(y_0) \subseteq L^+\bfP^1$ is dense open. At the generic point $\xi \in \Spec V$, we have
\begin{equation}\label{eq:generic_pt_of_V}
[f(\xi):1] = [1:f(\xi)^{-1}] : \xi \to D(x_0),
\end{equation}
with $f(\xi)^{-1} = \varpi ([t] + [a_0]\varpi + [a_1]\varpi^2 + [a_2] \varpi^3 + \dots)^{-1} = \sum_{n=1}^\infty [b_n]\varpi^n \in \widetilde K$ with some $b_n \in K$. Thus the image of \eqref{eq:generic_pt_of_V} lies in the closed subset $Z \subsetneq D(x_0)$ defined above.
In particular, $f(\xi) \neq f(s)=\eta$, and the proposition is proven. \qedhere
\end{proof}

\begin{rem}
In the equal characteristic case, including the case $\charac(k)=0$, Maffei obtained a different proof of \Cref{prop:P1_nonrep}, see the forthcoming work \cite{Maffei_unpub}. We also note that our proof of \Cref{prop:P1_nonrep} also applies when $k$ is of characteristic $0$ and $\Perf_k$ is replaced by the category of reduced $k$-algebras.
\end{rem}

Here is a related fact, which seems worth mentioning, but was not used in the above proof.

\begin{lm}\label{lm:ind_then_rep}
Let $X/F$ be a proper scheme. If $LX$ is strictly ind-representable, then it is representable.
\end{lm}
\begin{proof}
Let $LX = \dirlim_i X_i$ be a strict ind-presentation. Choose any proper model $\fX$ of $X$ over $\caO_F$. As the $k$-scheme $L^+\fX$ is qcqs, the natural map $L^+\fX \to LX$ factors through $X_n$ for some $n \in I$. As by \eqref{eq:1} $L^+\fX(\ff) \stackrel{\sim}{\to} LX(\ff) = \dirlim X_i(\ff)$ is a bijection for any algebraically closed $\ff/k$ and $X_i \to X_j$ are closed immersions for $j>i$, it follows that $X_n \to X_i$ is bijective on geometric points for $i\geq n$. But all $X_i$ are perfect, hence reduced, schemes and all transition maps are closed immersions, so $X_n \to X_i$ is an isomorphism for $i\geq n$. With other words, $LX = X_n$ is a scheme.
\end{proof}

\begin{cor}\label{cor:non_ind_rep_in_general}
Let $X$ be a scheme separated and of finite type over $F$. Suppose there exists a non-constant map $\bfP^1_{\breve F} \to X_{\breve F}$. Then $LX$ is not weakly ind-representable.
\end{cor}
\begin{proof}
As $LX_{\breve F} = LX \times_k \bar k$ (see \cite[1.a.1(a)]{PappasR_08}) we may replace $F$ by $\breve F$. Let then $g \colon \bfP^1 \to X$ be non-constant. As $X$ is separated over $F$, $g$ is proper. Also, $g$ is quasi-finite, as an infinite fibre would force it to be constant. Thus $g$ is finite. By the next lemma, $Lg$ is representable in affine morphisms, so that if $LX$ were a weak ind-scheme, then $L\bfP^1$ were too, contradiction.
\end{proof}

\begin{lm}
Let $g \colon Y \to X$ be a finite morphism of $F$-schemes. Then $Lg$ is representable in affine morphisms.
\end{lm}
\begin{proof}
Let $R \in \Perf_k$ and let $x: T:=\Spec R \to LX$ be any map corresponding to $\tilde x : \Spec \widetilde R\to X$.
As $g$ is finite, $Y\times_X \Spec \widetilde R = \Spec B$ for some finite $\widetilde R$-algebra $B$. We need to show that $T\times_{LX}LY$ is an affine $R$-scheme. For $R' \in \Perf_R$ compute
\begin{align*}
(T\times_{LX}LY)(R')&=\{y \in LY(R') : Lg(y) = x|_{R'}\} \\
&= \{\tilde y \in Y(\widetilde{R'}) \colon g(\tilde y) = \tilde x|_{\widetilde{R'}}\}\\
&= {\rm Hom}_{\Spec \widetilde R}(\Spec \widetilde{R'},Y\times_X \Spec \widetilde R)\\
&= {\rm Hom}_{\widetilde R-{\rm Alg}}(B,\widetilde{R'})
\end{align*}
Let $b_1,\dots,b_n$ be algebra generators of $B$ over $\widetilde R$. There is some $N>0$ such that $\varpi^Nb_i$ is integral over $\bW(R)$. Let $B^+$ be the $\bW(R)$-subalgebra of $B$ generated by $\varpi^N b_i$. Then $B^+$ is integral over $\bW(R)$ and $B^+[\varpi^{-1}]=B$.

Note that $\bW(R')$ is integrally closed within $\widetilde{R'}$. Indeed, let $x \in \widetilde{R'}$ be integral over $\bW(R')$ and write $x=\varpi^{-n}y$ with $y \in \bW(R')$ and $n\geq 0$ minimal. If $n>0$, multiply a monic equation of $x$ over $\bW(R')$  (which is, say, of degree $d>0$) by $\varpi^{nd}$, and reduce the result modulo $\varpi$. We get the equation $y^d \equiv 0$ in $\bW(R')/\varpi\bW(R') \cong R'$. But $R'$ is perfect, hence reduced, so $y\equiv 0 \,(\varpi)$, contradicting the minimality of $n$.

Inverting $\varpi$ induces a map ${\rm Hom}_{\bW(R)-{\rm Alg}}(B^+, \bW(R')) \to {\rm Hom}_{\widetilde R-{\rm Alg}}(B, \widetilde{R'})$. As $\bW(R') \subseteq \widetilde{R'}$, it is injective. We claim it is also surjective. Let $\alpha \in {\rm Hom}_{\widetilde R-{\rm Alg}}(B, \widetilde{R'})$. As $B^+$ is integral over $\bW(R)$, any element of $\alpha(B^+)$ will be integral over $\bW(R)$, and hence over $\bW(R')$. Thus by the above observation $\alpha(B^+) \subseteq \bW(R')$. With other words, $\alpha$ comes from ${\rm Hom}_{\bW(R)-{\rm Alg}}(B^+, \bW(R'))$.

It thus remains to show that $R' \mapsto {\rm Hom}_{\widetilde R-{\rm Alg}}(B,\widetilde{R'}) = {\rm Hom}_{\bW(R)-{\rm Alg}}(B^+, \bW(R'))$ is representable by $\Spec C$ for some $C\in \Perf_R$. But this functor is nothing else than $L_R^+(\Spec B^+)$, where for a $\bW(R)$-scheme $\fY$, $L^+_R\fY$ is the functor $R' \mapsto \fY(\bW(R'))$ on $\Perf_R$. Noting that $L_R^+(\Spec B^+) = L_k^+(\Spec B^+) \times_{L^+_k(\Spec \bW(R))} \Spec R$ (where the map $\Spec R \to L^+_k(\Spec \bW(R))$ corresponds to the identity of $\bW(R)$), it suffices to show that $L_k^+(\Spec B^+)$ is representable. By a result of Bhatt--Gabber, see \cite[Proposition 18.2]{BertapelleG_18_arxiv}, $L_k^+(\Spec B^+)$ agrees with the (perfect) Greenberg realization of the $\caO$-scheme $\Spec B^+$ in the terminology of \emph{loc.cit.} Then, by Proposition 18.5(iv) of \emph{loc.cit.}, $L^+(\Spec B^+)$ is then represented by an affine scheme.\qedhere
\end{proof}

\begin{rem} In the equal characteristic case the same non-representability statements hold also without perfection: if the loop space (regarded as a functor on all $k$-algebras) would be weakly ind-representable, then it would be so after perfection, contradicting the assertions proven above.
\end{rem}

In \cite[Theorem A]{MaffeiMV_26} Maffei--Melani--Vezzosi prove an ind-representability statement for five versions of two-dimensional affine Grassmannians for a solvable group $G$ over a two-local field of equal characteristic. They mention that their methods do not apply to general reductive groups. In fact, as a consequence of \Cref{prop:P1_nonrep}, for reductive groups the ind-representability result of \emph{loc.cit.} is already optimal for the \emph{loop of the usual affine Grassmannian} ${}^L\mathcal{G}r_{G}=L(\mathcal{G}r_G)$ (see Definition 1.5 of \emph{loc.cit.}):

\begin{cor}\label{cor:MaffeiMV_indrep}
Let $G$ be a reductive group over $k$. In the setup of \cite[Definition 1.5]{MaffeiMV_26} (in particular, $k=\bar{k}$, ${\rm char}(F)>0$), ${}^L\mathcal{G}r_G$ is ind-representable if and only if $G$ is a torus.
\end{cor}
\begin{proof}
The 'if' direction is proven in \emph{loc.cit.}
For the converse, suppose the semisimple rank of $G$ is at least one. Then there exists a closed immersion $\bfP^1 \har \mathcal{G}r_G$. As $L(\cdot)$ preserves closed immersions, $L\bfP^1 \har L(\mathcal{G}r_G)$ is also one. If $L(\mathcal{G}r_G)$ were ind-representable, $L\bfP^1$ were too. This would still be the case after perfection, contradicting \Cref{prop:P1_nonrep}.
\end{proof}

\section{Absence of rational curves implies representability}

The following is our main result.

\begin{thm}\label{thm:Lplus_equals_L}
Let $\fX$ be a proper $\caO$-scheme with generic fiber $X$.
Suppose that any $\bar k$-morphism $\bfP^1_{\bar k} \to \fX_{\bar k}$ is constant. Then the natural map $L^+\fX \stackrel{\sim}{\longrightarrow} LX$ is an isomorphism.
\end{thm}

\begin{rem}\label{rem:Johannes}
Note that \Cref{thm:Lplus_equals_L} is similar to \cite[Theorem 1.1]{Anschuetz_22}. Namely, when $\fX$ is either as in \Cref{thm:Lplus_equals_L} or $B\cG$ for a parahoric group scheme $\cG$ over $\caO_F$, then both results (the latter, and the evaluation of the first on complete valuation rings $V \in \Perf_k$ of rank one) just say that the natural map $\Hom(\caY_V, \fX) \to \Hom(\caY_V \sm \{s\}, \fX)$ is an equivalence of groupoids.
\end{rem}

In the course of the proof we will use the following lemmas.

\begin{lm}\label{lm:arc_bij_surj}
Let $f : X\to Y$ be a map of (set-valued) arc-sheaves on $\Perf_k$, and assume that $X$ is qcqs and $Y$ quasi-separated. Suppose that for each valuation ring $V \in \Perf_k$ of rank $\leq 1$ with algebraically closed fraction field, the map $f(V) \colon X(V) \to Y(V)$ is surjective (resp. bijective). Then $f$ is surjective (resp. an isomorphism).
\end{lm}
\begin{proof}
Throughout the proof 'surjective' means 'surjective as a map  of arc-sheaves'. We first prove the surjectivity assertion. It suffices to prove that for any $A \in \Perf_k$ and any map $T=\Spec A \to Y$, the map $T\times_Y X \to T$ is surjective. As $Y$ is quasi-separated and $T,X$ are quasi-compact, $T \times_Y X$ is quasi-compact. Thus there is a quasi-compact affine $Z \in \Perf_k$ and a surjection $Z\to T\times_Y X$. It suffices to show that $Z \to T$ is surjective, that is (as both are affine schemes) an arc-cover in the sense of \cite[Definition 1.2]{BhattM_18}. Let $V \in \Perf_k$ be a valuation ring of rank $\leq 1$ (we may wlog assume that $V$ is absolutely integrally closed). Given a map $a:\Spec V\to T$, its composition with $T \to Y$ lifts to a map $\Spec V \to X$ by the assumption of the lemma. This defines a lift $a'\colon \Spec V\to T\times_Y X$ of $a$. As $Z \to T\times_Y X$ is surjective, we may find some arc-cover $U \to \Spec V$ and a lift $b \colon U \to Z$ of $a'|_U$. By definition of an arc-cover applied to ${\rm id}_{\Spec V}$ we may (refining $U$) assume that $U=\Spec W$ for a valuation ring $W$ of rank $\leq 1$, which is faithfully flat over $V$. This shows that $Z \to T$ is an arc-cover, finishing the proof of the surjectivity claim. For the bijectivity claim, the proof of \cite[Lemma 2.1]{Ivanov_arc} applies.
\end{proof}

\begin{rem} We observe a mistake in \cite[Lemma 2.1]{Ivanov_arc}, which however does not affect any of the results of the article. For the surjectivity assertion, the direction (i) $\Rar$ (ii) claimed in \emph{loc.cit.} is trivially false: if $L/\bar{k}$ is an algebraically closed transcendental extension, then $\Spec L \to \Spec \bar k$ is a v-surjection of qcqs v-sheaves, but for the rank $0$ valuation ring $V=\bar k$, the map $(\Spec L)(V) = {\rm Hom}_{\bar k}(L,\bar k) = \varnothing \to \{{\rm id}\} = {\rm Hom}_{\bar k}(\bar k,\bar k) = (\Spec \bar k)(V)$ is not surjective. (Also the proof of the surjectivity assertion in (ii) $\Rar$ (i) of \emph{loc.cit.} has to be minimally corrected: cf. the proof of \Cref{lm:arc_bij_surj}.)
\end{rem}

\begin{lm}\label{lm:general_factoriz_modification}
Let $S= \Spec A$ be a local normal scheme with closed point $s$ and $S\sm\{s\} \neq \varnothing$. Let $T$ be an integral scheme and $q:T \to S$ a proper morphism, which is an isomorphism over $S \sm \{s\}$. Let $g: T \to Y$ be a morphism into a scheme. Suppose that $g(q^{-1}(s))$ is contained in an affine open subset of $Y$.
Then $g$ factors uniquely through $S$.
\end{lm}
\begin{proof}
Let $W = \Spec B \subseteq Y$ be an open affine neighborhood of $g(q^{-1}(s))$. The closed subset $T \sm g^{-1}(W)$ does not meet $q^{-1}(s)$, and hence $s\not\in q(T \sm g^{-1}(W))$. By properness of $q$, $q(T \sm g^{-1}(W)) \subseteq S$ is closed. As $S$ is local with closed point $s$, we necessarily must have $q(T \sm g^{-1}(W)) =  \varnothing$ and therefore $T \sm g^{-1}(W) =  \varnothing$. Thus $T \subseteq g^{-1}(W)$, and hence $g$ factors through $W$. By Stein factorization \cite[Tag 0AY8]{StacksProject}, one has $q_\ast \caO_T = \caO_S$. Thus the map $T \to W = \Spec B$ defines
\[B \to H^0(T,\caO_T) = H^0(S,q_\ast\caO_T) = A, \]
which gives the desired map $S \to W \subseteq Y$. Uniqueness also follows from $A=H^0(T,\caO_T)$.\qedhere
\end{proof}

\begin{proof}[Proof of \Cref{thm:Lplus_equals_L}]
$L^+\fX = \prolim_i L^+_i\fX$ is a cofiltered limit of qcqs schemes (see \cite{BertapelleG_18_arxiv}) along affine transition maps, hence itself representable by a qcqs scheme. Thus by \cite[Theorem 5.17]{BhattM_18} it is a qcqs arc-sheaf. By \Cref{prop:loop_arc} $LX$ is an arc-sheaf, too. By \cite[Lemma 3.5]{Ivanov_arc} $LX$ is separated. Thus by \Cref{lm:arc_bij_surj} it suffices to prove that for any valuation ring $V \in \Perf_k$ of rank $\leq 1$ with algebraically closed fraction field, $L^+\fX(V) \to LX(V)$ is bijective. By \eqref{eq:1} this holds when ${\rm rk} (V) = 0$, so we may assume ${\rm rk}(V) = 1$.
Choose a pseudo-uniformizer $t \in \fm_V \sm\{0\}$. The $t$-adic completion $V \to V^{\wedge}_t$ is an arc-cover. As $L^+\fX \to LX$ is injective by \Cref{lm:injective}, the bijectivity of $L^+\fX(V^{\wedge}_t) \to LX(V^{\wedge}_t)$ implies the bijectivity of $L^+\fX(V) \to LX(V)$. Thus we may assume that $V$ is $t$-adically complete.

Let $K$ (resp. $\kappa$) denote the fraction (resp. the residue field) of $V$. Let $x \in LX(V) = X(\widetilde V) = \fX(\widetilde V)$. Now $\Spec \widetilde V = D(\varpi) \subseteq \caY_V = \Spec \bW(V)$ is the locus where $\varpi$ is invertible. So we have to show that $x: D(\varpi) \to \fX$ extends uniquely to a morphism $\caY_V \to \fX$.

Observe that $\bW(V)_{(\varpi)}$ is the valuation subring of ${\rm Frac}\,\bW(V)$ for the $\varpi$-adic valuation, and hence a discrete valuation ring.
Thus by the valuative criterion of properness, the restriction of $x$ to ${D(\varpi) \cap \Spec \bW(V)_{(\varpi)}} $ extends uniquely to a morphism $x_1 : \Spec \bW(V)_{(\varpi)} \to \fX$. Choose an affine open neighborhood $\Spec B \subseteq \fX$ of the image of the closed point of $\Spec \bW(V)_{(\varpi)}$. Then $x_1$ factors through $\Spec B$, so we get a map $B \to \bW(V)_{(\varpi)}$ of $\caO$-algebras. As $\fX$ is of finite presentation over $\caO$, $B$ is too, so that map factors through $B \to \bW(V)_a$ for some $a$. With other words, $x_1$ extends to a map (again denoted $x_1$) on the open neighborhood $(\varpi) \in D(a) \subseteq \caY_V$. Moreover, $x_1$ agrees with $x$ on the overlap $D(a\varpi)$ (the locus where they agree is closed as $\fX$ is separated and contains the generic point of $\caY_V$ by construction). Thus $x$ extends (again, uniquely) to a morphism
\[x' : U:= \caY_V \sm \{s\} \subseteq  D(\varpi) \cup D(a) \to \fX,\]
where $s$ denotes the closed point of $\caY_V$.

Consider the graph $\gamma : U \to \caY_V \times_{\caO_F} \fX$ of $x'$. Let $\Gamma \subseteq \caY_V \times_{\caO_F} \fX$ be the scheme-theoretic image of $\gamma$. Denote the natural projections by $\caY_V \stackrel{p}{\leftarrow} \Gamma \stackrel{g}{\to} \fX$. By \cite[Tag 0CNG]{StacksProject} $p\colon \Gamma \to \caY_V$ is an isomorphism over $U$ (noting that $U \subseteq \caY_V$ is retrocompact).
As $\caY_V$ is reduced, \cite[056B]{StacksProject} shows that $\Gamma$ is too. Moreover, $\Gamma$ is integral as $U$ is. Finally, note that $p$ is proper, as $\fX/\caO$ is.
The key statement is the following.

\begin{lm}\label{lm:rational}
Any positive-dimensional irreducible component $C \subseteq \Gamma_{s,{\rm red}}$ of the reduced fibre of $\Gamma$ over $s$ is a rational (possibly singular) curve over $\kappa$.
\end{lm}
\begin{proof}
Let $\eta_C$ denote the generic point of $C$.
Here is a summary of the proof: We construct (step 1) a  rank one valuation $\Spec W \to \Gamma$ centered at $\eta_C$. Then we show (step 2) that the composed map $x:\Spec W \to \caY_V$ actually defines a rank-one point $x$ of the adic space $\caY_{\caS}$ from \cite[II.1.1]{FarguesScholze}, where $\caS = \Spa(K,V)$ is perfectoid over $k$. The residue field $k(\eta_C)$ embeds into the residue field $\widetilde{\caH(x)}$ of the completed residue field $\caH(x)$ at $x$ (step 3). Then we need to pass (step 4) to the standard perfectoid cover $\caY_{\caS} \times_{\Spa \caO} \Spa \caO_{\infty} \to \caY_{\caS}$ as in \cite[Proposition II.1.1]{FarguesScholze}\footnote{In \cite{FarguesScholze}, the setup only allows finite residue field $k$. However, \cite[II.1.1 Proposition]{FarguesScholze} holds for arbitrary $k$: the use of Frobenius in its proof to reduce to $n=1$ can be omitted by redoing the same calcuation for any $n$.},
and find a preimage $x_\infty$ of $x$ with tilt $x_\infty^\flat \in (\caY_{\caS} \times_{\Spa \caO} \Spa \caO_{\infty})^\flat \cong \bD_{\caS,{\rm perf}}$, where $\bD_{\caS,{\rm perf}}$ is the perfectoid open unit disc (this step can be omitted if $\charac(F) = p$). Finally (step 5), exploiting the preservation of residue fields under tilting, and the well-known computation of the residue fields at all four types of rank one points of the adic affine line over $\caS$, we can conclude.

For step 1, recall that $\Gamma \to \caY_V$ is birational, and let $\bW(V) \subseteq W_0 \subseteq {\rm Frac}({\bW(V)})$ be a valuation ring dominating the local ring $\caO_{\Gamma,\eta_C}$  (it exists by \cite[Tag 00IA]{StacksProject}). First of all, by properness of $\Gamma \to \caY_V$, we have a map $\Spec W_0 \to \Gamma$ sending the closed point to $\eta_C$. Write $I=(\varpi,[t]) \subseteq \bW(V)$. Since $W_0$ is a valuation ring, $IW_0$ is generated by a single element $a \in \{\varpi,[t]\}$. Let $\fq = \sqrt{aW_0}$ and $\fp = \bigcap_{i\geq 0} a^iW_0$. Then $W = (W_0/\fp)_{\fq/\fp}$ is a valuation ring of rank one (cf. \cite[Remark 2.2]{BhattM_18}) and we have the composed map
\[ \alpha \colon \Spec W \har \Spec W_0 \to \Gamma. \] The closed point $\bar{\fq} \in \Spec W \subseteq \Spec W_0$ specializes to the closed point of $\Spec W_0$, which in turn maps to $\eta_C \in \Gamma$. Thus $\bar{\fq}$ maps to a generization of $\eta_C$ in $\Gamma$. But the preimage of $\fq$ in $\bW(V)$ contains $I$ and hence also $\sqrt{I}$; thus $\bar\fq$ lies over $s \in \caY_V$, and therefore $\bar\fq$ maps into $\Gamma_s = p^{-1}(s)$ under $\Spec W_0 \to \Gamma$. As $\eta_C$ is a generic point of an irreducible component of $E$, we get $\alpha(\bar\fq)=\eta_C$. In particular we have a local map $\caO_{\Gamma,\eta_C} \to W$. Furthermore, replacing $W$ by its $a$-adic completion, we may assume that $W$ is complete.

For step 2, we have $IW = aW \neq 0$, thus either $\varpi$ or $[t]$ does not vanish in $W$. We claim that both do not vanish. Suppose first $\varpi \mapsto 0 \in W$. This gives a factorization $\Spec W \to \Spec V \to \caY_V$. As the generic point of $\Spec V$ lies in $U$, properness of $p$ gives a lift $\beta : \Spec V\to \Gamma$, whose composition with $\Spec W \to \Spec V$ agrees with $\alpha$ as $p$ is separated. Thus $\beta(s) = \eta_C$, which implies that the residue field of $\eta_C \in \Gamma$ equals $\kappa$, a contradiction as $\dim(C) \geq 1$. Thus $\varpi$ maps to a pseudo-uniformizer in $W$. Suppose next that $[t] \mapsto 0 \in W$. Then $[\fm_V]$ maps to $0$; as any element of $\ker(\bW(V) \to \bW(\kappa))$ has all Teichm\"uller coefficients in $\fm_V$, it maps into $\varpi^n W \subseteq W$ for each $n$, and as $\varpi \in W$ is a pseudo-uniformizer, we have $\bigcap_{n\geq 0} \varpi^n W = 0$. With other words, $\Spec W \to \caY_V$ factors through the closed subscheme $\Spec \bW(\kappa)$. Now $\bW(\kappa)$ is a valuation ring, and the same argument involving properness as above (with $\bW(\kappa)$ replacing $V$) gives a contradiction. Thus both $\varpi$ and $[t]$ map to non-zero in $W$. Let $x$ be the rank one valuation on $\bW(V)$ defined by the map $\bW(V) \to W$. This valuation is clearly continuous for the $(\varpi,[t])$-adic topology on $\bW(V)$: all topologically nilpotent elements of $\bW(V)$ lie in $\sqrt{I}$, hence map into $\bar\fq$, which consists of topologically nilpotent elements. We have shown that $x$ defines a point of the adic space $\caY_{\caS} = \Spa(\bW(V)) \sm V([t])$, where $\caS= \Spa(K,V)$.

For step 3, let $\fp_x \in \caY_V$ denote the support of $x$, so that $\caH(x)$ is the completion in ${\rm Frac}(W)$ of $k(\fp_x) = {\rm Frac}(\bW(V)/\fp_x)$. Let also $\xi \in \Gamma$ denote the image under $\alpha$ of the generic point of $\Spec W$; it lies over $\fp_x$. As $[t]$ maps to non-zero in $W$ (see step 2), $[t] \not\in \fp_x$ and hence $\Gamma \to \caY_V$ is an isomorphism around $\fp_x$. Thus $k(\xi) = k(\fp_x)$. As $\xi \rightsquigarrow \eta_C$, we have the map $\caO_{\Gamma,\eta_C} \to \caO_{\Gamma,\xi} \tar k(\xi)$. Its composition with $k(\xi) \har {\rm Frac}(W)$ coincides with our original map $\caO_{\Gamma,\eta_C} \to W \har {\rm Frac}(W)$. Thus its image is contained in the valuation ring $k(\xi) \cap W$ induced by $W$. As $O_{\Gamma,\eta_C} \to W$ is local, $\caO_{\Gamma,\eta_C} \to k(\xi) \cap W$ is. Thus it induces a $\kappa$-embedding of $k(\eta_C)$ into the residue field of $k(\xi) \cap W$. But the latter residue field coincides with the residue field $\widetilde{\caH(x)}$ of the completion $\caH(x)$ of $k(\xi)$. (If $\charac(F) >0$, we can now directly conclude by step 5 below.)

For step 4, denote by $F_\infty$ the completion of $F(\varpi^{\frac{1}{p^{\infty}}})$, and by $\caO_\infty$ its integers. Choose a continuous $F$-embedding $F_\infty \har \Omega=\widehat{\overline{\caH(x)}}$ into a completed algebraic closure of $\caH(x)$ (recall that by step 2, $x \not\in V(\varpi)$). Let $x_\infty \in \caY_{\caS} \times_{\caO}\caO_{\infty}$ be the image of the resulting $\Spa(\Omega,\Omega^\circ)$-point over $x$. Tilting via \cite[Proposition II.1.1]{FarguesScholze} gives $x_\infty^\flat \in (\caY_{\Spa(K,V)} \times_{\Spa \caO} \Spa \caO_{\infty})^\flat \cong \bD_{\caS,{\rm perf}}$, where $\bD_{\caS,{\rm perf}}$ is the perfectoid open unit disc over $\caS$, and
\[ k(\eta_C) \subseteq \widetilde{\caH(x)} \subseteq \widetilde{\caH(x_\infty)} \cong \widetilde{\caH(x_\infty^\flat)},\]
as tilting preserves residue fields by \cite[Lemma 3.4(iii)]{Scholze_PS}.

For step 5, let $\bD_{\caS}$ denote the usual adic open unit disc over $\caS$. The map $\bD_{\caS,{\rm perf}} \to \bD_\caS$ sends $x_\infty^\flat$ to a rank one point $y \in \bD_\caS \subseteq \bA^{1,{\rm ad}}_{\caS}$, so it is not a type-(5)-point in the usual classification of points of the adic affine line. Now $\caH(x_\infty^\flat)$ is the completion of the perfection of $\caH(y)$. As completion does not change the residue fields, $\widetilde{\caH(x_\infty^\flat)}=\widetilde{\caH(y)}^{\rm perf }$. In particular, $y$ cannot be of type (1), (3) or (4) either, as in this case the residue field $\widetilde{\caH(y)}$ would be algebraic over (and hence equal to) $\kappa$. Thus $y$ must be of type (2) -- a Gauss point of some disc -- and in this case it is well-known that $\widetilde{\caH(y)} \cong \kappa(z)$ is rational over $\kappa$. As perfection commutes with formation of residue fields, it follows that $\widetilde{\caH(x_\infty)} \cong \kappa(z^{\frac{1}{p^{\infty}}})$. Thus there exist embeddings $\kappa \subseteq k(\eta_C) \subseteq \kappa(z^{\frac{1}{p^{\infty}}})$. But $k(\eta_C)$ is finitely generated over $\kappa$, hence we even have $k(\eta_C) \subseteq \kappa(z^{\frac{1}{p^N}})$ for $N\gg 0$. By L\"uroth's theorem it follows that $k(\eta_C)$ is rational over $\kappa$. \qedhere
\end{proof}

By \cite[Tag 0AY8]{StacksProject}, $\Gamma_s = p^{-1}(s)$ is geometrically connected. By \Cref{lm:rational}, any irreducible component $C$ of $\Gamma_{s,{\rm red}}$ is either a point or a rational (possibly singular) curve over $\kappa$. Then, by assumption of the theorem, $\Gamma_s$ contracts to a single point of $\fX_k$ under $\Gamma \stackrel{g}{\to} \fX$. By \Cref{lm:normality}, $\caY_V$ is normal. Combining all these facts with the properties of $\Gamma,U$ and $p$ stated before \Cref{lm:rational}, we see that \Cref{lm:general_factoriz_modification} is applicable and shows that $\Gamma \to \fX$ factors uniquely through a desired map $\caY_V \to \fX$.\qedhere
\end{proof}

The idea of the following remark (as well as the main idea of the proof of \Cref{thm:Lplus_equals_L}) is very similar to \cite[Proposition 6.2]{GabberLL_15} (the analogous situation for a discrete valuation ring $V$ is essentially a special case of \emph{loc.cit.}).

\begin{rem}\label{rem:Temkin_eq_char_case}
In the equal characteristic case, \Cref{thm:Lplus_equals_L} has a purely schematic proof, which we quickly sketch here. Namely, the $U$-modification $\Gamma \to \caY_V$ from the proof can be refined by a $U$-admissible blowup $T\to \Gamma \to \caY_V = \Spec V[\![\varpi]\!]$ by \cite[Tag 081T]{StacksProject}. Being admissible, it algebrizes to an admissible blowup $C' \to C=\bfA^1_V$ which is an isomorphism outside the closed point $x=\sqrt{(\varpi,t)}\in C$ (where $t \in V$ pseudo-uniformizer). Now, by a result of Temkin \cite[Proposition 4.5.1]{Temkin_sst} applied to the curve $C' \to \Spec V$ (with empty multisection $D$), there is a stable modification $C^{\rm st} \to C'$, which is (by Corollary 1.3(ii) of \emph{loc.cit.}) an isomorphism over the semistable locus. Thus the composition $C^{\rm st} \to C' \to C$ is an isomorphism outside $x$. Now, by Lemma 4.4.2 of \emph{loc.cit.} applied to $C^{\rm st} \to C$, the fiber over $x$ in $C^{\rm st}$ is a $\bfP^1$-tree or a point; in particular, it is connected and each positive-dimensional irreducible component is a $\bfP^1$. Now base changing along the flat (by \Cref{lm:flatness_rk1val}) map $\caY_V \to C$, and using that blowups commute with flat base change, $C^{\rm st} \times_C \caY_V \to T \to \caY_V$ is a $U$-modification whose fiber over $x$ is either a point or a connected tree of $\bfP^1$'s. Now one can conclude by \Cref{lm:general_factoriz_modification} as in the proof of the theorem.
\end{rem}

\subsection{Changing the model}\label{sec:changing_model}
Let $\caC_\caO$ denote the category of proper $\caO$-schemes with no rational curves in the geometric special fibre. As a consequence of \Cref{thm:Lplus_equals_L} we deduce that the functor $\fX \mapsto L^+\fX$ on $\caC_\caO$ factors through the generic fibre.

\begin{cor}\label{cor:isogeny}
 The jet scheme functor \[\fX \mapsto L^+\fX \,\, \colon \,\,\caC_\caO \to \{\text{perfect $k$-schemes}\}\] factors through the generic fibre functor $\fX \mapsto \fX \times_\caO F$ from $\caC_\caO$ to the full subcategory of all proper $F$-schemes which admit a model $\fX \in \caC_\caO$.
 \end{cor}
\begin{proof}
This follows directly from \Cref{thm:Lplus_equals_L}.
\end{proof}

\begin{rem}\label{rem:conseq}
In particular, if $\fX, \fX' \in \caC_\caO$ are two proper $\caO$-models of a proper scheme $X/F$, then $L^+\fX \simeq L^+\fX'$ are canonically isomorphic $k$-schemes. However, their truncations might differ: in general one has $L^+_r \fX \not\simeq L^+_r\fX'$.
\end{rem}

\begin{rem}\label{rem:Loeser_Sebag}
It would be interesting to study the relation of \Cref{thm:Lplus_equals_L} and \Cref{cor:isogeny} with the motivic integration theory on rigid-analytic varieties over $F$ of Loeser--Sebag \cite{LoeserS_03}. In particular, in \cite[Theorem-Definition 4.1.2]{LoeserS_03} they prove independence of the motivic integral from the choice of the model $\fX$. Now (a formal scheme analogue of) \Cref{cor:isogeny} would show that the space over which the integral is taken already does not depend on the model (at least if its geometric special fiber satisfies the assumption of \Cref{thm:Lplus_equals_L}). However, giving a new approach to this result would require bringing the differential form $\omega$ into the picture.
\end{rem}

\bibliography{bib_ADLV}{}

\providecommand{\bysame}{\leavevmode\hbox to3em{\hrulefill}\thinspace}
\providecommand{\MR}{\relax\ifhmode\unskip\space\fi MR }
\providecommand{\MRhref}[2]{%
  \href{http://www.ams.org/mathscinet-getitem?mr=#1}{#2}
}
\providecommand{\href}[2]{#2}
\begin{thebibliography}{MMV26}

\bibitem[AL25]{AnschuetzLB_25_Fourier}
Johannes Ansch{\"u}tz and Arthur-C{\'e}sar {Le Bras}, \emph{A {Fourier}
  transform for {Banach--Colmez} spaces}, J. Eur. Math. Soc. \textbf{27}
  (2025), no.~9, 3651--3712.

\bibitem[And21]{Andreychev_21}
Grigory Andreychev, \emph{Pseudocoherent and perfect complexes and vector
  bundles on analytic adic spaces}, preprint (2021),
  \href{https://arxiv.org/abs/2105.12591}{arXiv:2105.12591}.

\bibitem[Ans22]{Anschuetz_22}
Johannes Ansch{\"u}tz, \emph{Extending torsors on the punctured
  {$\operatorname{Spec}(A_{\mathrm{inf}})$}}, J. reine und angew. Math.
  \textbf{783} (2022), 227--268.

\bibitem[BGA22]{BertapelleG_18_arxiv}
Alessandra Bertapelle and Cristian~D. Gonz{\'a}lez-Avil{\'e}s, \emph{{The
  Greenberg Functor Revisited}}, 2022, Version 6,
  \href{https://arxiv.org/abs/1311.0051v6}{arXiv:1311.0051v6}.

\bibitem[Bha16]{Bhatt_16}
Bhargav Bhatt, \emph{Algebraization and {Tannaka} duality}, Camb. J. Math.
  \textbf{4} (2016), no.~4, 403--461.

\bibitem[BM21]{BhattM_18}
Bhargav Bhatt and Akhil Mathew, \emph{The arc-topology}, Duke Math. J.
  \textbf{170} (2021), no.~9, 1899--1988.

\bibitem[Bou06]{Bourbaki_CommAlg}
Nicolas Bourbaki, \emph{{Alg{\`e}bre commutative. Chapitres 5 {\`a} 7}},
  {\'E}l{\'e}ments de math{\'e}matique, Springer, Berlin, 2006.

\bibitem[Dri03]{Drinfeld_03}
V.~Drinfel'd, \emph{Infinite-dimensional vector bundles in algebraic geometry
  (an introduction)}, Preprint (2003),
  \href{https://arxiv.org/abs/math/0309155}{arXiv:math/0309155}.

\bibitem[FS26]{FarguesScholze}
Laurent Fargues and Peter Scholze, \emph{Geometrization of the local
  {L}anglands correspondence}, Ast{\'e}risque, no. 466, Soc. Math. France,
  Paris, 2026.

\bibitem[GLL15]{GabberLL_15}
Ofer Gabber, Qing Liu, and Dino Lorenzini, \emph{Hypersurfaces in projective
  schemes and a moving lemma}, Duke Math. J. \textbf{164} (2015), no.~7,
  1187--1270.

\bibitem[GLQ23]{GrosSQ_23}
Michel Gros, Bernard {Le Stum}, and Adolfo Quir{\'o}s, \emph{{Cartier}
  transform and prismatic crystals}, Tunisian J. Math. \textbf{5} (2023),
  no.~3, 405--456.

\bibitem[Iva23]{Ivanov_arc}
Alexander~B. Ivanov, \emph{{Arc-descent for the perfect loop functor and
  $p$-adic Deligne--Lusztig spaces}}, J. reine angew. Math. (Crelle)
  \textbf{794} (2023), 1--54.

\bibitem[KV04]{KapranovV_04}
Mikhail Kapranov and Eric Vasserot, \emph{Vertex algebras and the formal loop
  space}, Publ. Math. IH{\'E}S \textbf{100} (2004), 209--269.

\bibitem[LS03]{LoeserS_03}
Fran{\c{c}}ois Loeser and Julien Sebag, \emph{Motivic integration on smooth
  rigid varieties and invariants of degenerations}, Duke Math. J. \textbf{119}
  (2003), no.~2, 315--344.

\bibitem[Lur17]{Lurie_HA}
Jacob Lurie, \emph{Higher algebra}, 2017,
  \href{https://www.math.ias.edu/~lurie/papers/HA.pdf}{https://www.math.ias.edu/~lurie/papers/HA.pdf}.

\bibitem[Lur26]{Lurie_Kerodon}
\bysame, \emph{Kerodon}, \url{https://kerodon.net}, 2026.

\bibitem[Maf]{Maffei_unpub}
Andrea Maffei, \emph{{The space of polynomial arcs in $\mathbb{P}^1$ is not an
  ind-scheme}}, forthcoming.

\bibitem[Mat89]{Matsumura_CRT}
Hideyuki Matsumura, \emph{Commutative ring theory}, 2nd ed., Camb. Studies in
  Adv. Math., vol.~8, Cambridge University Press, Cambridge, 1989.

\bibitem[Mat22]{Mathew_FaithfullyFlatDescent_22}
Akhil Mathew, \emph{{Faithfully flat descent of almost perfect complexes in
  rigid geometry}}, J. Pure Appl. Algebra \textbf{226} (2022), no.~5, 106938.

\bibitem[MMV26]{MaffeiMV_26}
Andrea Maffei, Valerio Melani, and Gabriele Vezzosi, \emph{{Two dimensional
  versions of the affine Grassmannian and their geometric description}}, Forum
  of Mathematics, Sigma \textbf{14} (2026), e42.

\bibitem[PR08]{PappasR_08}
George Pappas and Michael Rapoport, \emph{Twisted loop groups and their affine
  flag varieties}, Adv. Math. \textbf{219} (2008), 118--198.

\bibitem[Ric20]{Richarz_20_basics}
Timo Richarz, \emph{{Basics on Affine Grassmannians}}, Lecture notes,
  \href{https://timo-richarz.com/wp-content/uploads/2020/02/BoAG\_02.pdf}{https://timo-richarz.com/wp-content/uploads/2020/02/BoAG\_02.pdf},
  2020.

\bibitem[Sch12]{Scholze_PS}
Peter Scholze, \emph{Perfectoid spaces}, Publ. Math. IH\'ES \textbf{116}
  (2012), no.~1, 245--313.

\bibitem[Sch26a]{Scholze_Berkovich_motives}
\bysame, \emph{Berkovich motives}, J. Amer. Math. Soc. \textbf{39} (2026),
  no.~3, 697--764.

\bibitem[Sch26b]{Scholze_ECD}
Peter Scholze, \emph{{\'E}tale cohomology of diamonds}, Ast{\'e}risque, no.
  466, Soci{\'e}t{\'e} Math{\'e}matique de France, Paris, 2026, In
  \emph{Geometrization of the local Langlands correspondence}, pp.~439--659.

\bibitem[Shi15]{Shimomoto_15}
Kazuma Shimomoto, \emph{On the {Witt} vectors of perfect rings in positive
  characteristic}, Communications in Algebra \textbf{43} (2015), no.~12,
  5328--5342.

\bibitem[Shi17]{Shimomoto_17}
\bysame, \emph{An embedding problem of noetherian rings into the {Witt}
  vectors}, preprint (2017),
  \href{https://arxiv.org/abs/1503.02018}{arXiv:1503.02018}.

\bibitem[{Sta}26]{StacksProject}
The {Stacks Project Authors}, \emph{{\itshape {S}tacks {P}roject}},
  \url{http://stacks.math.columbia.edu}, 2026.

\bibitem[SW20]{ScholzeW_20}
Peter Scholze and Jared Weinstein, \emph{{Berkeley} lectures on {$p$}-adic
  geometry}, Ann. Math. Studies, vol. 207, Princeton University Press,
  Princeton, NJ, 2020.

\bibitem[Tem10]{Temkin_sst}
Michael Temkin, \emph{Stable modification of relative curves}, J. Algebraic
  Geom. \textbf{19} (2010), 603--677.

\bibitem[Zhu17]{Zhu_17}
Xinwen Zhu, \emph{Affine {G}rassmannians and the geometric {S}atake in mixed
  characteristic}, Ann. Math. \textbf{185} (2017), no.~2, 403--492.

\end{thebibliography}
\bibliographystyle{amsalpha}
\end{document}